\documentclass[12pt,fleqn]{article}
\usepackage{graphicx}
\usepackage{mathptmx}
\usepackage{amsmath,amssymb,amsfonts,amsthm}
\usepackage{latexsym}
\def\pmb#1{\setbox0=\hbox{$#1$}%
\kern-.025em\copy0\kern-\wd0
\kern.05em\copy0\kern-\wd0
\kern-.025em\raise.0433em\box0}

\newfont{\lie}{eufm10 at 12pt}
\newfont{\field}{msbm10 at 11pt}
\newcommand{\R}{\mathbb{R}}
\newcommand{\h}{\mbox{\lie h}}                  % \h = gotic h

\newcommand{\sm}[1]{\mbox{\small $#1$}}

\newcommand{\La}[1]{\mbox{\Large $#1$}}
\newcommand{\LA}[1]{\mbox{\LARGE $#1$}}

\newtheorem{Theorem}{Theorem}
\newtheorem{Lemma}[Theorem]{Lemma}
\newtheorem{Corollary}[Theorem]{Corollary}
\newtheorem{Proposition}[Theorem]{Proposition}
\begin{document}
\setcounter{page}{1}

\markright{\sl \hfill Li -- Salavessa \hfill}
\title{ \normalsize \bf  GRAPHIC BERNSTEIN RESULTS IN
CURVED PSEUDO-RIEMANNIAN MANIFOLDS}
\author{Guanghan Li$^{1,2\dag}$ and Isabel M.C.\ Salavessa$^{2\ddag}$}
\date{}
\protect\footnotetext{\!\!\!\!\!\!\!\!\!\!\!\!\! {\bf MSC 2000:}
Primary: 53C21, 53C42, 53C50
$\dag$ Partially supported by NSFC (No.10501011), NSF of Hubei 
(No. 2008CDB009), and by
Funda\c{c}\~{a}o Ci\^{e}ncia e Tecnologia (FCT) through a FCT
fellowship SFRH/BPD/26554/2006.\\
$\ddag$ Partially supported by FCT through the Plurianual of CFIF
and POCI-PPCDT/MAT/60671/2004.}
\maketitle ~~~\\[-5mm]
{
\footnotesize $^1$ School of Mathematics and Computer Science,
Hubei University, Wuhan, 430062, P. R. China,
e-mail: liguanghan@163.com}\\[2mm]
{\footnotesize $^2$ Centro de F\'{\i}sica das Interac\c{c}\~{o}es
Fundamentais, Instituto Superior T\'{e}cnico, Technical University
of Lisbon, Edif\'{\i}cio Ci\^{e}ncia, Piso 3, Av.\ Rovisco Pais,
1049-001 Lisboa, Portugal;~~
e-mail: isabel.salavessa@ist.utl.pt}\\[5mm]
\begin{abstract}
\noindent
We generalize  a  Bernstein-type result
due to Albujer and Al\'ias,
 for  maximal surfaces in a
curved Lorentzian
product 3-manifold of the form $\Sigma_1\times \mathbb{R}$,
to higher dimension and codimension.
We consider $M$  a complete spacelike graphic submanifold
with parallel mean curvature,
defined by a map $f: \Sigma _1\rightarrow \Sigma _2$
between two Riemannian
manifolds $(\Sigma _1^m, g_1)$ and $(\Sigma^n _2, g_2)$
of sectional curvatures $K_1$ and $K_2$, respectively.
We take on $\Sigma_1\times \Sigma _2$ the
pseudo-Riemannian product metric $g_1-g_2$. Under the curvature
conditions,  $\mathrm{Ricci}_1\geq 0$ and
$K_1\geq  K_2$, we prove  that,
if the second fundamental form of $M$  satisfies an integrability condition,
then $M$ is totally geodesic, and it is a slice if
$\mathrm{Ricci}_1(p)>0$ at some point. For bounded
 $K_1$,  $K_2$ and  hyperbolic angle $\theta$,
 we conclude $M$ must be maximal. If $M$ is a maximal
 surface and $K_1\geq K_2^+$,  we show $M$
is totally geodesic with no need for
further assumptions.
Furthermore, $M$ is a slice if at some
point $p\in \Sigma_1$, $K_1(p)> 0$, and if  $\Sigma_1$ is flat
and $K_2<0$ at some point $f(p)$, then the
image of $f$ lies on a geodesic  of $\Sigma_2$.\\[3mm]
\end{abstract}
\section{Introduction and statement of main results}

The classical Bernstein theorem states that an entire minimal graph in
$\R^3$ is a plane. This result has been generalized to  graphic
hypersurfaces of $\R^{m+1}$ for
$m\le 7$, and for higher dimensions and codimensions under various
growth conditions.  Calabi in \cite {c} introduced  a similar problem in
 Minkowski space. He considered a maximal (that is, with mean curvature
$H=0$) spacelike hypersurface $M$ in the Lorentz-Minkowski space
$\R_1^{m+1}$ with
metric
$ds^2=\sum _{i=1}^m(dx_i)^2-(dx_{m+1})^2$,
given by the graph of a function
$f$ on $\R^m$ with $|Df|<1$. In this case the equation for maximality takes
the form
\[
\sum _{i=1}^m\frac {\partial }{\partial x_i}
\left(\frac {{\partial f}/{\partial
x_i}}{\sqrt{1-|Df|^2}}\right)=0.
\]
Calabi, for $m\leq 4$, and later Cheng and Yau  \cite {cy1} for any $m\geq 2$,
proved that any entire solution to the
above equation is linear.

Spacelike hypersurfaces of constant mean curvature  in Lorentzian spaces
have been the subject of investigation in general relativity theory
(see for example \cite{ars}).
Contrarily to the Riemannian case, there are entire spacelike graphs with
nonzero constant mean curvature in $\mathbb{R}^{m+1}_1$ (see  \cite{tr}),
 as for example the hyperboloids.
On the other hand,
under some  boundedness assumptions or  growth conditions
 on  the Gauss map, 
 Bernstein-type results have been obtained
for spacelike submanifolds 
 of the pseudo-Euclidean space $\R^{m+n}_n$ with parallel
mean curvature (\cite {x2,xy,jx2}).

A natural generalization is to consider maximal
spacelike  submanifolds in a non-flat ambient
space. Albujer and Al\'{\i}as  \cite {aa}  proved a
new Calabi-Bernstein-type result for surfaces immersed into a
Lorentzian product 3-manifold of the form $\Sigma_1\times
\mathbb{R}$, where $\Sigma_1$ is a Riemannian surface
of nonnegative Gauss curvature.
In this note, we generalize this result to
  spacelike graphic  submanifolds 
 with parallel mean curvature in a non-flat pseudo-Riemannian
product space of any dimension $m+n$, and under less restrictive
curvature conditions.
Our main tool, as in \cite{aa},  is the explicit computation 
 of  $\Delta\cosh\theta$, where $\theta$ is the hyperbolic angle, 
a quantity
 that measures how far is
$M$ from a slice, and that was introduced
by Chern \cite{ch2}  for the Riemannian version of the Bernstein
theorem for  surfaces in $\mathbb{R}^3$.
 This angle plays a similar role
to the Gauss map in flat ambient spaces mentioned above.
The main difficulty in higher codimension is the 
definition of the hyperbolic angle itself.  This can be done with
the help of a suitable  parallel form 
from the ambient space (see section 2), 
representing a different approach comparing with
the case $m=2$ and $n=1$, namely on what concerns the
 computations involving $\cosh\theta$, that are also considerably
 more complicated when $n\geq 2$, even for $m=2$.

Let $N=\Sigma_1\times \Sigma _2$ be a pseudo-Riemannian product
manifold of two Riemannian manifolds $(\Sigma_i, g_i)$ with
pseudo-Riemannian metric $\bar{g}=g_1-g_2$. We denote by
$K_i$ and $\mathrm{Ricci}_i$ the  sectional curvatures and
the Ricci tensor of each $\Sigma_i$, respectively.
Assume $\Sigma_1$ is oriented
and $M$ is a spacelike graphic submanifold
$\Gamma_{f}=\{(p,f(p)): p\in \Sigma_1\}$, defined by a  map
$f:\Sigma_1\rightarrow\Sigma _2$ with $f^*g_2<g_1$.
Let $g$ be the induced metric on $M$.
The \em hyperbolic angle \em $\theta$
can be defined by
\begin{equation}\label{(1)}
\cosh\theta=\frac{1}{\sqrt{\det(g_1-f^*g_2)} },
\end{equation}
where the determinant is taken  with respect to $g_1$.
Note that $\|df\|^2$ is bounded.
Let $B$ be the second
fundamental form of $M$.
 We  state our first main theorem:
\begin{Theorem} Assume $M$ is  a complete spacelike graph
with parallel mean curvature, and for any $p\in \Sigma_1$, 
 $\mathrm{Ricci}_1(p)\geq 0$ and
  $K_1(p)\ge K_2(f(p)) $.  Then we have:\\[2mm]
$(i)$~  If  $\|df\|\, \|B\|$ is an integrable function
on $M$ (and this is the case $M$ compact),
 then $M$ is totally geodesic.
Moreover, if
$\mathrm{Ricci}_1(p)>0$ at some point, then $M$
is a slice, that is, $f$ is constant.\\[1mm]
$(ii)$~  If  $K_1$,
$K_2$ and  $\cosh\theta$ are bounded, then $M$ is maximal. Furthermore, 
if $M$ is non-compact, and for $p$ away from a compact set,
 $K_1(p)-K_2(f(p))\geq \varepsilon$ ($\mathrm{Ricci}_1(p)\geq \varepsilon$, 
respectively), 
where $\varepsilon >0$ is a constant,
 then $f$ cannot have rank greater or equal to
two (one, respectively) at infinity.
\end{Theorem}
\noindent
In case $\Sigma_2$ is
one-dimensional, in next proposition
 we replace  the
boundedness  condition on $\cosh\theta$ and $K_1$ given in  theorem 1$(ii)$  
by  weaker conditions, and obtain the same result,
by recalling an isoperimetric inequality in \cite{sal2}:
\begin{Proposition}
If $\Sigma_2$ is one dimensional,  $\Sigma_1$ is complete,
$\mathrm{Ricci}_1\geq 0$, and $\cosh\theta= o(r)$ when
$r\rightarrow +\infty$,
where $r$ is the distance function to a point of $\Sigma_1$,
then $M$ is maximal.
\end{Proposition}
For $M$ a  Riemannian surface, using parabolicity arguments
we obtain:
\begin{Theorem} If $M$ is a complete maximal spacelike graphic surface, 
and for each  $p\in \Sigma _1$,
$K_1(p)\ge\max\{0, K_2(f(p))\} $,  then $M$ is totally geodesic.
 Furthermore:\\[1mm]
$(i)$~~~ If $K_1(p)>0$ at some point $p\in M$,  then $M$ is a slice. \\
$(ii)$~~  (\cite{c} if $n=1$, \cite{jx2})
If $\Sigma_1=\mathbb{R}^2$ and $\Sigma_2=\mathbb{R}^n$, then $M$ is a plane.
\\
$(iii)$~ If $\Sigma_1$ is flat
and $K_2<0$ at some point $f(p)$, then either
$M$ is a slice or the image of $f$ is a geodesic of $\Sigma_2$.
\end{Theorem}
\noindent
Our proof in $(ii)$ gives partially a simpler proof of the same
result of Jost and Xin in \cite{jx2} for the case of surfaces.
As a consequence of theorem 3, if  $\Sigma _2=\mathbb{R}$ we have:
\begin{Corollary}[\cite{aa}] Let $M$ be a complete maximal spacelike surface of
$N=\Sigma _1\times \mathbb{R}$, with pseudo-Riemannian product
metric $g_1-dt^2$, and assume $M$ can be written as the graph of a smooth
map $f:\Sigma_1\rightarrow \mathbb{R}$. If
$K_1\geq 0$  then $M$ is totally geodesic.
Moreover, if $K_1>0$ at some point of $\Sigma_1$, then $M$ is a
slice.
\end{Corollary}
This paper is organized as follows. In section 2, we recall some
preliminaries of spacelike submanifolds in pseudo-Riemannian manifolds,
and  compute the Laplacian of
$\cosh\theta $. The proofs of  theorem 1 and proposition 2 are given in
section 3. In section 4, we discuss the surface case and prove
Theorem 3.

\section{Spacelike submanifolds in pseudo-Riemannian products}

Let $N$ be a $(m+n)$-dimensional pseudo-Riemannian manifold with
non-degenerate metric  $\bar{g}$  of index $n$, and  
$F: M\rightarrow N$  a $m$-dimensional
spacelike submanifold immersed into $N$. We
denote by $\overline{\nabla} $, $\nabla$ and ${\nabla}^{\bot}$
 the  connections on
$N$, $M$, and the normal bundle $NM$, respectively,
 and by $B$  the second fundamental
form of $M$.
We convention the sign of the curvature tensor
$\bar{R}$ of $N$ is defined by
$\bar{R}(X, Y)=[\overline{\nabla}_X,\overline{\nabla}_Y]
-\overline{\nabla}_{[X,Y]}$
and
$\bar{R}(X, Y,Z, W)=\bar{g}(\bar{R}(Z, W)Y, X)$.
 We make use of the indices
range, $i, j, k, \cdots, =1, 2, \cdots, m$, $\alpha, \beta,
\cdots, =m+1, \cdots, m+n$, and $a, b, c, \cdots, =1, 2, \cdots,
m+n$, and 
 choose  orthonormal frame fields $\{e_1, \cdots, e_{m+n}\}$
of $N$, such that restricting to $M$, $\{e_1, \cdots,
e_m\}$ is a tangent frame  (of spacelike vectors), and $\{e_{m+1}, \cdots,
e_{m+n}\}$ a normal frame (of timelike vectors).
Let $h_{ij}^{\alpha}$ be the components of the 
second fundamental form, $B(e_i,
e_j)=h_{ij}^{\alpha}e_{\alpha}$, and $\bar{R}^a_{bcd}$,
 $R^i_{jkl}$, $R^{\alpha}_{\beta kl}$,
 the components of
the curvature tensors of $N$, $M$ and $NM$, respectively.
Using the  structure equations of
$N$, we derive (see e.g.\ \cite{cy1,jx2} for the case $N=\mathbb{R}^{n+m}$)
\begin{eqnarray}
&&R_{jkl}^i=\bar{R}_{jkl}^i-\sum
_{\alpha}(h_{ik}^{\alpha}h_{jl}^{\alpha}-h_{il}^{\alpha}h_{jk}^{\alpha})
~~~~(\mbox{Gauss~equation})\label{(2.2)}\\
&&R_{\beta kl}^{\alpha}=\bar{R}_{\beta
kl}^{\alpha}-\sum_i(h_{ki}^{\alpha}h_{li}^{\beta}-h_{li}^{\alpha}
h_{ki}^{\beta}),~~~(\mbox{Ricci~equation})\nonumber\\
&&h_{ij,k}^{\alpha}-h_{ik,j}^{\alpha}
=-\bar{R}_{ijk}^{\alpha}=\bar{R}_{\alpha ijk}. ~~~~~~~~~~
(\mbox{Codazzi~equation})\label{(2.4)}
 \end{eqnarray}
The mean curvature of $F$ is denoted by
$H =trace_g B= H^{\alpha}e_{\alpha}$, $H^{\alpha}=\sum_{i}h^{\alpha}_{ii}$.
Let $\Omega$ be a parallel $m$-form on $N$.  Similarly to \cite{wa2},
 we  compute the Laplacian of the pull back $F^*\Omega$,
$\Delta F^*\Omega=\sum_k \nabla_k\nabla _kF^*\Omega-
\nabla_{\nabla_{e_k}e_k}F^*\Omega$.
\begin{eqnarray} \label{(3.1)}
(\nabla _kF^*\Omega)(e_{1}, \cdots, e_{m})&=&\sum_i\Omega(e_1,\ldots,
(\overline{\nabla} _ke_{i}-\nabla _ke_{i}), \ldots,e_{m} ) \nonumber\\
&=&\sum_i\Omega(e_1,\ldots, B(e_k,e_i),\ldots, e_m).
\end{eqnarray}
Set $u=u^{\top}+u^{\bot}$, with $u^{\top}\in TM$ and $u^{\bot}\in NM$.
Differentiating (\ref{(3.1)}), we have
\begin{eqnarray*}
\Delta F^*\Omega(e_{1}, \cdots, e_{m}) =&&\!\!
\sum_{ik}\Omega(e_1,\ldots,\!\nabla_{e_k}^{\bot}B(e_k,e_i)
\!+\!(\overline{\nabla}_{e_k}B(e_k,e_i))^{\top}\!\!,\ldots,e_m)\\
&&+\sum_k\sum_{j<i}\Omega( e_1, \ldots,
B(e_k,e_j),\ldots,B(e_k,e_i),\ldots,e_m)\\
&&+\sum_k\sum_{j>i}\Omega( e_1, \ldots,
B(e_k,e_i),\ldots,B(e_k,e_j),\ldots,e_m).\\
\end{eqnarray*}
Using  (\ref{(2.4)}),
$\sum_k\nabla^{\bot}_{e_k}B(e_k,e_i)=
\nabla^{\bot}_{e_i}H+(\bar{R}(e_k,e_i)e_k)^{\bot},$
and that
$$\sum_{ik}g((\overline{\nabla}{e_k}B(e_k,e_i))^{\top},e_i)=
\sum_{ik}-\bar{g}(B(e_k,e_i),B(e_k,e_i))= \|B\|^2,$$
 we get in components
\begin{equation}
(\Delta F^*\Omega )_{1\cdots m}=\Omega _{1\cdots
m}||B||^2+2\sum_{\alpha <\beta, i<j}\!\!\!\!\Omega _{\alpha \beta ij}
\hat{R}_{\beta ij}^{\alpha }
+\sum _{\alpha, i}\Omega _{\alpha i}H_{,i}^{\alpha }-\sum
_{\alpha, i, k}\Omega _{\alpha i}\bar{R}_{kik}^{\alpha },~~~
\label{(3.2)}
\end{equation}
where $\hat{R}_{\beta ij}^{\alpha}=\sum_kh_{ik}^{\alpha
}h_{jk}^{\beta}-h_{ik}^{\beta}h_{jk}^{\alpha}$, and $\Omega
_{\alpha \beta ij}=\Omega (e_1, \cdots , e_{\alpha },\cdots,
e_{\beta},\cdots,e_m)$ with $e_{\alpha}$, $e_{\beta}$ occupying
the $i$-th and the $j$-th positions. The same meaning is for
$\Omega _{\alpha i}$. Here, $||\cdot||$  denotes the
absolute of the norm of a timelike vector.

Next we consider the pseudo-Riemannian
product manifold  $N=\Sigma_1\times \Sigma
_2$, with the pseudo-Riemannian
 metric $\bar{g}=g_1-g_2$,  where
 $(\Sigma_i, g_i)$ are two Riemannian manifolds
of dimension $m$ and $n$, respectively. Let $\pi_i:N\to \Sigma_i$
denote the corresponding projections.

Suppose $M$ is a spacelike graph of a map $f: \Sigma
_1\rightarrow \Sigma _2$.
For any $p\in \Sigma _1$, we consider $\lambda_1^2\geq \lambda_2^2\geq
\ldots\geq \lambda_m^2\geq 0$ the eigenvalues of
$f^*g_2$.
The spacelike condition on $M$
means $\lambda_i^2<1$.
By the  Weyl's perturbation theorem \cite{we},
ordering the eigenvalues in this way, each
$\lambda_i^2:\Sigma_1\rightarrow [0,1)$ is a continuous
locally Lipschitz function.
 For each $p$, 
let $ s=s(p)\in\{1\,\ldots, m\} $ be the rank of $f$ at $p$, that is,
$\lambda_s^2>0$ and $\lambda_{s+1}^2=\ldots=\lambda_m^2=0$.
Then $s\leq \min\{m,n\}$.
We say that \em  $f$ has rank $\geq s$ at infinity, \em  if
there is a constant $\epsilon>0$ such that  $\lambda_s^2\geq \epsilon$,
away from a compact set $K$. 
We take an  orthonormal basis
$\{a_i\}_{i=1, \cdots, m}$ of
$T_p\Sigma _1$ of eigenvectors of $f^*g_2$ with
corresponding eigenvalues $\lambda_i^2$.
Set $a_{i+m}=df(a_i)/|df(a_i)|$
for $i\leq s$. This constitutes an orthonormal system in
$T_{f(p)}\Sigma_2$, that we complete to give
an orthonormal basis
 $\{a_{\alpha}\}_{\alpha =m+1, \cdots, m+n}$ for
$T_{f(p)}\Sigma _2$.
Moreover, changing signs of the $\lambda_i$ if necessary, we can write
 $df(a_i)=-\lambda _{i\alpha}a_{\alpha}$,
where $\lambda _{i\alpha }=\delta_{\alpha,m+i}\lambda_i$
 meaning $=0$ if $i>s$
or $\alpha> m+s$.
Therefore
\begin{eqnarray} \label{(3.4)}
e_i &=&\sm{\frac 1{\sqrt{1-\sum
_{\beta}\lambda _{i\beta}^2}}}(a_i+\sum _{\beta}\lambda
_{i\beta}a_{\beta}),~~~~{i=1, \cdots , m} \\
\label{(3.5)}
e_{\alpha}&=&\sm{\frac 1{\sqrt{1-\sum _j\lambda
_{j\alpha }^2}}}(a_{\alpha}+\sum _j\lambda
_{j\alpha}a_{j}),~~~\alpha =m+1, \cdots , m+n
\end{eqnarray}
form orthonormal basis of $T_pM$ and  $N_pM$, respectively. We may assume
$e_i$ to be positively
oriented.
We also identify $M=\Gamma_f$ with $\Sigma_1$ with the graph metric
$g=g_1-f^*g_2$, and
 consider $\lambda_i$ as functions
on the variable $x=(p,f(p))\in M$ identified with the variable
$p\in \Sigma_1$,  through the diffeomorphism
${\pi_1}_{|_M}:M\rightarrow \Sigma_1$.
Let $\Omega $ be the volume form of
$(\Sigma _1, g_1)$, which is a parallel $m$-form on $N$,
and $\Omega'$ the one of $(\Sigma_1,g)$.
The ratio  $\Omega/\Omega'$  is given by
\[
\Omega _{1\cdots m}= *({\pi_1}_{|_M})^*\Omega= \pi_1^*\Omega(e_1,\ldots, e_m)=
\frac 1{\sqrt{\prod _{i=1}^{m}(1-\lambda
_i^2)}}=\frac{1}{\sqrt{det(g_1-f^*g_2)}},
\]
where $*$ is the star operator on $M$.
This quantity is $\geq 1$ and is  $\cosh\theta$ 
defined in (\ref{(1)}).
Then,  $\cosh\theta$ is identically equal to
$1$ if and only if $f$ is a constant map, that is, $M$ is a slice.
Furthermore,  if $\cosh\theta$ is bounded, what means
$\lambda_i^2\leq 1-\delta$, where $\delta> 0$ is a constant,
then 
$g_1\geq g\geq \delta g_1$. In this case, 
 $\Sigma_1$ is complete (compact) if and only if  $M$ is so.
The   singular values $\lambda_i$
of $f$ are constant maps if $F$ is a totally geodesic immersion.
To see this, we first note that $F$ is  totally geodesic
if and only if
 $f:\Sigma_1\rightarrow \Sigma_2$ is a totally geodesic map
(a proof of this is similar to the Riemannian case
\cite {sal1}, remark 2).
Parallel transport of $a_i$ along  geodesics
 of $\Sigma_1$
 starting
from $p$, shows that $f^*g_2(a_i,a_j)=\delta_{ij}\lambda_i^2$ is  constant.

With respect to the frames $\{e_i\}$  and
$\{e_{\alpha}\}$ given in (\ref{(3.4)}), (\ref{(3.5)}), we have
\begin{equation}\label{(3.6)}
2\sum
_{\alpha <\beta, i<j}\Omega _{\alpha \beta ij}\hat{R}_{\beta
ij}^{\alpha }=2\sum _{\alpha ,\beta , k, i<j }\lambda _{i\alpha
}\lambda
_{j\beta}(h_{ik}^{\alpha}h_{jk}^{\beta}-h_{ik}^{\beta}h_{jk}^{\alpha})
\cosh\theta.
\end{equation}
We denote by $R_1$ and $R_2$ the curvature tensors of $\Sigma _1$
and $\Sigma _2$.
Since
$\bar{R}_{kik}^{\alpha}=\sum_{\beta}\bar{R}_{\beta
kik}\bar{g}^{\alpha \beta}=-\bar{R}_{\alpha kik}$, we get
 \begin{eqnarray*}
\lefteqn{-\bar{R}_{kik}^{\alpha}=
\bar{R}(e_{\alpha},e_k, e_i, e_k)=}\\
&=& R_1(\pi_1(e_{\alpha}),
\pi_1(e_k),\pi_1(e_i),\pi_1(e_k))-R_2(\pi_2(e_{\alpha}),
\pi_2(e_k),\pi_2(e_i),\pi_2(e_k))\\
&=&{ \frac{\La{\{}{\sum
_{l}}\lambda_{l\alpha}R_1(a_l, a_k, a_i, a_k)-{\sum
_{\beta,\gamma,\delta }}\lambda_{k\beta}\lambda_{i\gamma}\lambda
_{k\delta} R_2(a_{\alpha},a_{\beta}, a_{\gamma},
a_{\delta})\La{\}}}
{\sqrt{(1-\sum_j\lambda _{j\alpha }^2)(1-\lambda
_{i}^2)}(1-\lambda _{k}^2)}}.
\end{eqnarray*}
Consider for $i\neq j$ the two-planes $P_{ij}=span\{a_i,a_j\}$,
$P'_{ij}=span\{a_{m+i},a_{m+j}\}$.
Since $\lambda_{i\alpha}$ is diagonal, we have
\begin{eqnarray}
\lefteqn{-\sum _{\alpha, i, k}\Omega _{\alpha
i}\bar{R}_{kik}^{\alpha}=}\nonumber\\
&=& \sum_{i,j}\sm{{\frac{\cosh\theta\, \lambda _i^2}{(1-\lambda
_i^2)(1-\lambda _j^2)}}} \left(R_1(a_i, a_j, a_i, a_j)-\lambda
_j^2R_2(a_{m+i},a_{m+j}, a_{m+i},a_{m+j})\right).~\nonumber\\
&=& \cosh\theta\sum_{i\neq j}\left({\sm{\frac{ \lambda _i^2}{(1-\lambda
_i^2)}}}K_1(P_{ij})+\sm{\frac{ \lambda _i^2\lambda_j^2}
{(1-\lambda_i^2)(1-\lambda _j^2)}} [K_1(P_{ij})-K_2(P'_{ij})]\right)
\label{(3.7)}
\end{eqnarray}
Inserting (\ref{(3.6)}) and (\ref{(3.7)}) into (\ref{(3.2)}),
and using the fact that the star operator is parallel,
we at last arrive at
\begin{eqnarray}
\lefteqn{\Delta \cosh\theta=
 \cosh\theta \LA{\{}||B||^2+2\sum _{k,i<j}\lambda
_i\lambda _jh_{ik}^{m+i}h_{jk}^{m+j}-2\sum _{k,i<j}\lambda
_i\lambda _jh_{ik}^{m+j}h_{jk}^{m+i}}\nonumber\\
&& +\sum_{i}\LA{(}{\sm{\frac{ \lambda _i^2}{(1-\lambda
_i^2)}}}\mathrm{Ricci}_1(a_i,a_i)+\sum_{j\neq i}
\sm{\frac{ \lambda _i^2\lambda_j^2}
{(1-\lambda_i^2)(1-\lambda _j^2)}} [K_1(P_{ij})-K_2(P'_{ij})]\LA{)}
\LA{\}}\nonumber\\
&&+\sum _{\alpha, i}\Omega _{\alpha i}H_{,i}^{\alpha},
\label{(3.9)}
\end{eqnarray}
\section{Proof of main results}
It is convenient
to recall lemma 3.1 in
\cite {aa}, which is also valid in our setting.
\begin{Lemma}[\cite{aa}, \rm{for} $n=2$, $m=3$]
Let $M$ be a spacelike  $m$-dimensional submanifold immersed into
$N=\Sigma _1\times \Sigma _2$. Then $\Sigma _1$ is necessarily
complete if $M$ is complete. In this case ${\pi_1}_{|_M}:
M\rightarrow \Sigma_1$
is a covering map.
\end{Lemma}
We also recall the well-known Omori-Cheng-Yau
maximum principle:
\begin{Proposition}[\cite {cy2}]
Let $u$ be a function bounded from above on a complete
manifold $M$ with Ricci curvature bounded from below. Then for any
fixed $p_0\in M$, there exists a sequence of points $\{p_k\}\subset
M$, such that $u(p_0)\le u(p_k)$, $~\lim _{k\rightarrow \infty}u(p_k)
=\sup_M u$, $\lim_{k\rightarrow \infty}|\nabla u(p_k)|=0$, and 
$\lim_{k\rightarrow \infty}\Delta u(p_k)\le 0.$
\end{Proposition}
\noindent
{\bf Proof of Theorem 1.}
By (\ref{(3.1)}) we have
\begin{equation}
\label{bound}
d\cosh\theta(e_k)=\sum_{i\alpha}\Omega (\pi(e_1),
\ldots,\pi_1(e_{\alpha}),  \ldots,
\pi_1(e_m))h_{ik}^{\alpha}
=\cosh\theta \sum _i\lambda _ih_{ik}^{m+i} 
\end{equation}
which implies
\begin{equation}\label{(4.1)}
\frac {|\nabla \cosh\theta|^2}{\cosh^2\theta
}=\sum _k(\sum _i\lambda _ih_{ik}^{m+i})^2=\sum
_{i,k}(\lambda_ih_{ik}^{m+i})^2+2\sum _{i<j,k}\lambda _i\lambda
_jh_{ik}^{m+i}h_{jk}^{m+j}.
\end{equation}
We shall calculate
$~
\Delta \ln(\cosh\theta)=(
{\cosh\theta\, \Delta(\cosh\theta)-|\nabla \cosh\theta |^2})
\, {\cosh^{-2}\theta.
}~$
From (\ref{(3.9)}) and (\ref{(4.1)}), and the assumption $H$ parallel,
that is $H_{,i}^{\alpha }=0$, we have, 
\begin{eqnarray}
\lefteqn{\Delta \ln(\cosh\theta)=||B||^2-\sum
_{i,k}\lambda_i^2(h_{ik}^{m+i})^2-2\sum _{k,i<j}\lambda
_i\lambda _jh_{ik}^{m+j}h_{jk}^{m+i}} \label{(4.3)}\\
&& +\sum_{i}\left(\sm{{\frac{ \lambda _i^2}{(1-\lambda
_i^2)}}}\mathrm{Ricci}_1(a_i,a_i)+\sum_{j\neq i}
\sm{\frac{ \lambda _i^2\lambda_j^2}
{(1-\lambda_i^2)(1-\lambda _j^2)}} [K_1(P_{ij})-K_2(P'_{ij})]\right)
\label{(4.4)}
\end{eqnarray}
First we need to compute the terms on the right hand side of (\ref{(4.3)}).
Since $F$ is spacelike, at each point $p\in M$
there exists a positive
 constant $\delta=\delta(p) \leq 1$ such that $\lambda _i^2\le 1-\delta$
 for any $1\le i\le m$. Thus, $|\lambda _i\lambda _j|\leq 1-\delta $
for any $i$ and $j$.
 We note that $\lambda _i=0$ for $i>\min (m,n)$. Therefore, we have
\[||B||^2\ge \sum _{i,k,j}(h_{ik}^{m+j})^2=\sum
_{i<j,k}[(h_{ik}^{m+j})^2+(h_{jk}^{m+i})^2]+\sum
_{i,k}(h_{ik}^{m+i})^2,\]
 where we keep in mind that
$h_{ik}^{m+j}=0$ when $m+j>m+n$ (because it is the only possible
meaning).
 So the terms in (\ref{(4.3)}) satisfy at $p$
\begin{eqnarray}
\quad&&||B||^2-\sum _{i,k}\lambda_i^2(h_{ik}^{m+i})^2-2\sum
_{k,i<j}\lambda _i\lambda
_jh_{ik}^{m+j}h_{jk}^{m+i}\nonumber\\
&\ge&\delta ||B||^2+(1-\delta)\left\{\sum
_{i<j,k}\left[(h_{ik}^{m+j})^2+(h_{jk}^{m+i})^2\right]+\sum _{i,
k}(h_{ik}^{m+i})^2\right\}\nonumber\\
&&-(1-\delta)\sum _{i,k}(h_{ik}^{m+i})^2-2(1-\delta)\sum
_{k,i<j}|h_{ik}^{m+j}||h_{jk}^{m+i}|\nonumber\\
&\ge& \delta ||B||^2.\label{(4.5)}
\end{eqnarray}
From the curvature assumptions,
$(\ref{(4.4)})\geq 0$.
Consequently,
 we have at each $p$ a
differential inequality,
\begin{equation}\label{(4.6)}
\Delta \ln
(\cosh\theta )\ge \delta(p) ||B||^2\ge \frac {\delta(p)}{m}||H||^2.
\end{equation}
Now we prove $(i)$.  From (\ref{(4.6)}), $\ln(\cosh\theta )$ is
a subharmonic function on a complete Riemannian manifold $M$,
and by (\ref{bound}), 
 $\|\nabla\ln(\cosh\theta )\|\leq C\|df\|\, \|B\|$,
where $C>0$ is a constant. Under the integrability condition of
$\|df\|\, \|B\|$, we have integrability of $\|\nabla\ln(\cosh\theta )\|$, 
and applying  Yau's special Stokes' theorem (corollary of \S1 \cite{yau})
we conclude that $\Delta \ln(\cosh\theta )=0$. From (\ref{(4.6)})
 we obtain  $||B||^2=0$, and therefore $M$ is totally geodesic.
Then all singular values $\lambda_i$ are constant functions.
Moreover if there exists at least one
point $p_0\in \Sigma _1$ such that $ \mathrm{Ricci}_1(p_0)>0$, then we easily
obtain from (\ref{(4.4)})  that $\lambda _i=0$ for any $i=1,
\cdots, m$, and thus $f$ is a constant map, that is, $M$ is a slice.

$(ii)$ Using an o.n. basis $E_i$
that at a given point diagonalizes the Ricci tensor $\mathrm{Ricci}^M$  of 
$M$, and applying Gauss equation (\ref{(2.2)}),  we have for each $s$
\begin{eqnarray}
\lefteqn{\mathrm{Ricci}^M(E_s,E_s)=\sum_{j\ne s}\{\bar{R}(E_s,E_j,E_s,E_j)
-\sum_{\alpha}(h_{ss}^{\alpha}h_{jj}^{\alpha}
-h_{sj}^{\alpha}h_{sj}^{\alpha})\}~~~}\nonumber\\
&=&\sum _{j\ne s}\bar{R}(E_s,E_j,E_s,E_j)+\sum
_{\alpha}\left(h_{ss}^{\alpha}-\frac 12{H^{\alpha
}}\right)^2
-\frac 14||H||^2+\sum _{\alpha , j\ne s}(h_{sj}^{\alpha})^2,~~~~
\label{(4.7)}
\end{eqnarray}
where all components appearing
in this expression are with respect to this frame.
Since $M$ has parallel mean curvature, $||H||$ is 
constant, and so $\mathrm{Ricci}^M$ is bounded from below
 whenever $\sum_{j}\bar{R}(E_s,E_j,E_s,E_j)$
is so. Using the other frame, we have
$\sum_j\bar{R}(E_s,E_j,E_s,E_j)
=\sum_{jik}  A_{si}A_{sk} \bar{R}(e_i,e_j,e_k,e_j)$, where
 $A_{si}=g(E_s,e_i)$, defines an orthogonal matrix. 
 As in section 2, we have
\begin{eqnarray}\label{(4.9)}
\lefteqn{\sum _{j\ne i,k}\bar{R}_{ijkj}=
\sum _{j\ne i, k}\bar{R}(e_i, e_j,e_k,e_j)=}\nonumber\\
&=& \sum _{j\ne i}{\frac{\La{(}R_1(a_i,a_j,a_k,a_j)
-\lambda_i \lambda_k\lambda
_j^2R_2(a_{m+i},a_{m+j},a_{m+k},a_{m+j})\La{)}}
{\sqrt{(1-\lambda _i^2)(1-\lambda_k^2)}(1-\lambda_j^2)}}.
\end{eqnarray}
Since $K_1$ and $K_2$ are bounded, the same holds for
 $R_1(a_i,a_j,a_k,a_j)$ and
$R_2(a_{m+i},$ $a_{m+j},a_{m+k},a_{m+j})$.
Boundedness of $\cosh\theta$ means there is a positive lower bound
$\delta$ for all $\delta(p)$,
and (\ref{(4.9)}) is bounded.  Thus,  from
Proposition 6,  there exists a sequence
$\{p_k\}\subset M$ such that
$\lim _{k\rightarrow \infty}\Delta \ln(\cosh\theta)(p_k) \le 0.$
Then by (\ref{(4.6)}) we conclude $H=0$, because $\|H\|$ is constant.
Now we assume for $p$ away from a compact set $K$ of $\Sigma_1$,
$K_1(p)-K_2(f(p))\geq \varepsilon$, where $\varepsilon>0$ is constant.
Again by (\ref{(4.5)}) and (\ref{(4.4)}), we have for $p\notin K$,
\begin{eqnarray}
\Delta \ln(\cosh\theta)\geq  \delta ||B||^2+ \sum_{i}
\La{(}\sm{\frac{\lambda_i^2}{(1-\lambda _i^2)}}\mathrm{Ricci}_1(a_i,a_i)
+\sum_{j\neq i}\sm{\frac{\lambda_i^2\lambda _j^2}{(1-\lambda _i^2)
(1-\lambda _j^2)}}\varepsilon\La{)}\ge 0.~~~~~
\label{(4.8)}
\end{eqnarray}
First we assume  $p_k\notin K$  for $k$ large.
From (\ref{(4.8)})
$\lim_{k\rightarrow \infty}\lambda _i(p_k)\lambda_j(p_k)=0$,
for $i\neq j$, and so
 $\lim_{k\rightarrow \infty}\lambda _i(p_k)=0$ for  $i\geq 2$.
In particular, $f$ cannot have rank greater or equal to two
at infinity.
If $p_k\in K$ for a subsequence, then $\ln(\cosh\theta)$
 attains a maximum at some limit point
$p_{\infty}$  of $p_k$,
and being a subharmonic map (see (\ref{(4.6)})),
by the strong maximum principle $\cosh\theta$ is constant,
and so, by (\ref{(4.6)}),  $B=0$, because $\delta(p_{\infty})\neq 0$. 
Consequently all $\lambda_i$ are  constant, and 
$0=(\ref{(4.8)})$ holds  for $p\notin K$, what implies
 $\lambda_i(p)=0$ for all $i\geq 2$. The case $\mathrm{Ricci}_1\geq
\varepsilon$ is similar.
\qed\\[4mm]
We observe the boundedness condition on $\theta$  is equivalent to
the boundedness condition on the Gauss map of Jost and Xin \cite{jx2}
in case $K_1=K_2=0$. Hence,  by their Bernstein theorem we have:
\begin{Proposition}[\cite{jx2}] If $\Sigma_1=\mathbb{R}^m$, $\Sigma_2=
\mathbb{R}^n$,  and $M$ is a graphic parallel
submanifold with bounded hyperbolic angle, then $M$ is a plane.
\end{Proposition}
Next we prove proposition 2. We are assuming $\Sigma_1$ complete and
oriented. Let $\mathcal{O}$ be an open set of $\Sigma_1$.
The  Cheeger constant of  $\mathcal{O}$ is defined by
\[\h (\mathcal{O})=\inf_D \frac{A_1(\partial D)}{V_1(D)},\]
where $D$ ranges over all open submanifolds of $\mathcal{O}$ with compact
closure in $\mathcal{O}$ and smooth boundary, and $A_1(\partial D)$ 
and $V_1(D)$
are respectively, the induced volumes of $\partial D$ and $D$ for the
metric  $g_1$. If $\Sigma_1$ is closed, we  adopt the same definition
for $\h(\Sigma_1)$, that is zero in this case. 
We fix a point $p$ of $\Sigma_1$,  
and $B_s$ denotes the open geodesic  ball at $p$ of radius $s$.
\begin{Lemma}[\cite{adr}]
If $\mathrm{Ricci_1}\geq 0$, then for any $r>0$,
$ \h(B_r)\leq \frac{C}{r}$,
where $C>0$ is a constant that does not depend on $s$.
In particular $\h(\Sigma_1)=0$.
\end{Lemma}
\begin{Proposition}[\cite{sal2}] If $\Sigma_2$ is one-dimensional
and  $M$ is a graphic spacelike hypersurface $\Gamma_f$,
then on a open bounded set $D$ of $\Sigma_1$, with smooth
boundary 
\[\inf_D\|H\|\leq \frac{1}{m}\frac{b_D}{\sqrt{1-b_D^2}}\frac{A_1(\partial D)}
{V_1(D)},~~~~\mbox{where}~~b_D=\sup_D\|\nabla f\|_1. \]
In particular, if $M$ has constant mean curvature, 
the Cheeger constant
of $(\Sigma_1,g_1)$ vanish,  and the hyperbolic angle is bounded, then
$M$ is maximal.
\end{Proposition}
Now, proposition 2 follows directly as an application of previous lemma and the
inequality in proposition 9, by taking $D\subset B_r$, and that for $n=1$,
 $\cosh\theta=(1-\|\nabla f\|^2)^{-1/2}$.

We should note that a 
 certain  nonnegativeness condition on the curvature of $\Sigma_1$
 plays a fundamental role in this type
of results.
If $\Sigma_1$ is the $m$-hyperbolic space $\mathbb{H}^m$ there are
examples of complete entire graphic hypersurfaces with constant mean
curvature $c$,  for  any $c$,
and with bounded hyperbolic angle,
as can be shown by the following proposition. The function
$r(x)=\ln \left(\frac{1+|x|}{1-|x|}\right)$
is the distance function
in $\mathbb{H}^m$ to $0$, for the Poincar\'{e} model:
\begin{Proposition}[\cite{sal2}]
Let $c$ be any constant and $f_c:\mathbb{H}^m\rightarrow \mathbb{R}$
defined by:
\[f_c(x)=\int_0^{r(x)}
\frac{\frac{c}{(\sinh r)^{m-1}}\int_0^r(\sinh t)^{m-1}dt}
{\sqrt{1+\left(\frac{c}{(\sinh r)^{m-1}}
\int_0^r(\sinh t)^{m-1}dt\right)^2}} dr. \]
Then $f_c$ is smooth on all
$\mathbb{H}^m$, and for each $c, d\in \mathbb{R}$, $\Gamma_{(f_c)+d}\subset
\mathbb{H}^m\times \mathbb{R}$ is a complete spacelike graph
of bounded hyperbolic angle,
with $~|\nabla f_c|^2_1\leq \frac{c^2}{(m-1)^2}/(1+\frac{c^2}{(m-1)^2})<1$~
 and  constant mean curvature given by
$~\langle H,\nu\rangle =\frac{c}{m}~$, where
$\nu=$ 
$(-\nabla f_c,1)/\sqrt{1+|\nabla f_c|^2_1} $ is the unit
 timelike normal to the graph. Furthermore,
$\{\Gamma_{(f_c)+d}(x):~ x\in \mathbb{H}^m, d\in \mathbb{R}\}$ (with $c$ fixed)
and $\{\Gamma_{(f_c)+d+c}(x):~ x\in \mathbb{H}^m, c\in \mathbb{R}\}$
(with $d$ fixed)
define foliations of $\mathbb{H}^m\times \mathbb{R}$ by hypersurfaces,
 with the same constant mean curvature $c$, and with constant mean curvature
parameterized by the leaf, respectively.
\end{Proposition}
\noindent
The constant mean curvature
parameter $c$ (or $1/c$) in the second example of the proposition,
 can be interpreted as
a natural ``time function''  of geometric nature. The existence
of such foliations have been considered  in  general relativity.
For $c=0$ the above examples are slices.
In  \cite{aa} there are several existence theorems and  explicit
examples of complete maximal graphic surfaces in
$\mathbb{H}^2\times \mathbb{R}$ that are not slices.
\section{Surface case $(m=2)$}
{\bf Proof of Theorem 3.}
We calculate
\[\Delta \left(\frac 1{\cosh\theta}\right)=
-\frac {\Delta \cosh\theta}{(\cosh\theta)^2}
+\frac{2|\nabla \cosh\theta|^2}{(\cosh\theta)^3}.\]
So, by (\ref{(3.9)}) and (\ref{(4.1)}) we have,
\begin{eqnarray}
\Delta\left( \frac 1{\cosh\theta}\right)
&=& -\frac 1{\cosh\theta}
\LA{( }||B||^2-2\sum _{k,i<j}\lambda _i\lambda
_jh_{ik}^{m+i}h_{jk}^{m+j} \nonumber\\
&&\quad  -2\sum _{k,i<j}\lambda
_i\lambda _jh_{ik}^{m+j}h_{jk}^{m+i}-2\sum _{i,k}\lambda _i^2(h_{ik}^{m+i})^2
\nonumber\\
&&\quad
+ \sum_{i=1,2}\sm{\frac{\lambda _i^2}{(1-\lambda _i^2)}}K_1
+\sm{\frac{\lambda_1^2\lambda_2^2}{(1-\lambda_1^2)(1-\lambda_2^2)}}
[K_1-K_2(a_3,a_4)]\LA{)}.
\label{(5.1)}
\end{eqnarray}
Since $M$ is  maximal  and  $m=2$, we have
\begin{eqnarray}
&&||B||^2-2\sum _{k,i<j}\lambda _i\lambda
_jh_{ik}^{m+i}h_{jk}^{m+j}-2\sum _{k,i<j}\lambda
_i\lambda _jh_{ik}^{m+j}h_{jk}^{m+i}-2\sum _{i,k}\lambda _i^2(h_{ik}^{m+i})^2
\nonumber\\
&=&||B||^2-2\lambda _1\lambda
_2\left[h_{11}^{m+1}h_{12}^{m+2}+h_{12}^{m+1}h_{22}^{m+2}+h_{11}^{m+2}
h_{12}^{m+1}+h_{12}^{m+2}h_{22}^{m+1}\right]\nonumber\\
&&-2\sum _{i,k}\lambda _i^2(h_{ik}^{m+i})^2\nonumber\\
&\ge&
\sum_{i<j,k}[(h_{ik}^{m+j})^2+(h_{jk}^{m+i})^2]+\sum_{i,k}(h_{ik}^{m+i})^2
-2\sum _{i,k}\lambda _i^2(h_{ik}^{m+i})^2\nonumber\\
&\ge& \sum_k (h_{1k}^{m+2})^2+(h_{2k}^{m+1})^2-(h_{1k}^{m+1})^2
-(h_{2k}^{m+2})^2
\nonumber\\
&=&0. \label{(5.2)}
\end{eqnarray}
Therefore by assumption on the curvatures, (\ref{(5.1)}) becomes
\begin{equation}\label{(5.3)}
\Delta\left( \frac 1{\cosh\theta}\right)\leq
-\frac 1{\cosh\theta}
\left(\sum_{i=1,2}\sm{\frac{\lambda _i^2}{(1-\lambda _i^2)}}K_1
+\sm{\frac{\lambda_1^2\lambda_2^2}{(1-\lambda_1^2)(1-\lambda_2^2)}}
[K_1-K_2(a_3,a_4)]\right)\le 0.
\end{equation}
By Gauss equation the
Gauss curvature of $M$ is given by
\[K_M=R_{1212}=\bar{R}_{1212}-\sum _{\alpha}(h_{11}^{\alpha}
h_{22}^{\alpha}-(h_{12}^{\alpha})^2)=\bar{R}_{1212}+
\sum_{\alpha}[{(h_{11}^{\alpha}})^2+(h_{12}^{\alpha})^2],\]
 where similarly to (\ref{(4.9)}),  
\begin{eqnarray*}
&&\bar{R}_{1212}=\frac 1{(1-\lambda _1^2)(1-\lambda
_2^2)}\left[K_1-\lambda_1^2\lambda _2^2K_2(a_3,a_4)\right]\geq 0.
\end{eqnarray*}
Consequently,  the Gauss
curvature of $M$ is nonnegative, and so $M$ is parabolic,
in the sense that any nonnegative superharmonic function on the surface
is constant.
By (\ref{(5.3)}), $\cosh\theta$ is constant, and the inequalities in
(\ref{(5.2)})
are identities. From these identities,
we immediately have
\[h_{ij}^{\alpha}=0 \quad {\mbox for} \quad \alpha \ge 5\] and
\[\sum_{i,k=1}^2(h_{ik}^{2+i})^2=\sum_{i,k=1}^2\lambda _i^2(h_{ik}^{2+i})^2.
\]
Since $\lambda _i<1$ for $i=1,2$, the last equality implies
$h_{ij}^3=h_{ij}^4=0$. Therefore, $M$ is totally geodesic and
 $\lambda_i$ and $\cosh\theta$ are constant.
From (\ref{(5.3)}) we conclude that, 
if at some point $K_1(p)>0$,  then $f$ is constant.
Now assume $K_1$ is identically  zero,
and $K_2<0$ at some point $f(p)$.
 By (\ref{(5.3)}) $\lambda_1\lambda_2= 0$.
Hence,  the  rank of $f$ is zero or one, and
 since $f:\Sigma_1\rightarrow
\Sigma_2$ is a totally geodesic map,
 the image of $f$ lies on a geodesic of $\Sigma_2$.
In case $\Sigma_i$ are  Euclidean spaces, a totally geodesic
surface of $\mathbb{R}^{2+n}_n$ is a plane.
\qed\\[4mm]

So we have the following corollary of Theorem 3
for $K_2=-1$:
\begin{Corollary} If $M$ is a complete maximal spacelike graph
of the pseudo--Riemannian product
$\mathbb{R}^2\times \mathbb{H}^n$, defined by
 a map $f:\mathbb{R}^2\rightarrow \mathbb{H}^n$, then
 either $f$ is constant or its image lies on a
geodesic of $\mathbb{H}^n$.
\end{Corollary}
\noindent
If $\Sigma_2$ is complete,
there are trivial  examples of complete totally geodesic spacelike graphs
in $\mathbb{R}^m\times \Sigma_2$,
with image of $f$ a non-constant geodesic.
Let  $\gamma:\mathbb{R}\rightarrow \Sigma_2$
be an entire geodesic
with $\|\gamma'(0)\|^2<1$. Then
$f:\mathbb{R}^m\rightarrow \Sigma_2$, given by
$f(x_1,\ldots,x_m)=\gamma(x_1)$, is a totally geodesic map with image
 $\gamma$, and the graph of $f$ is a complete totally
geodesic spacelike immersion.
Note that we are using two facts: geodesics of $\Sigma_2$
are just the  totally geodesic maps from $\mathbb{R}$ into $\Sigma_2$,
and totally geodesic graphs are just the graphs
of totally geodesic maps  (see section 2).

 \end{document}